# Maximin setting for investment problems and fixed income management with observable but non-predictable parameters

Nikolai Dokuchaev

Department of Mathematics and Computer Science,

The University of West Indies, Mona, Jamaica.

email ndokuch@uwimona.edu.jm

**Abstract**

We study optimal investment problem for a diffusion market consisting of a finite number of risky assets (for example, bonds, stocks and options). Risky assets evolution is described by Itô's equation, and the number of risky assets can be larger than the number of driving Brownian motions. We assume that the risk-free rate, the appreciation rates and the volatility of the stocks are all random; they are not necessary adapted to the driving Brownian motion, and their distributions are unknown, but they are supposed to be currently observable. Admissible strategies are based on current observations of the stock prices and the aforementioned parameters. The optimal investment problem is stated as a problem with a *maximin* performance criterion. This criterion is to ensure that a strategy is found such that the minimum of utility over all distributions of parameters is maximal. Then the *maximin* problem is solved for a very general case via solution of a linear parabolic equation.

## 1 Introduction

This paper studies optimal investment problem for a diffusion market consisting of a finite number of risky assets (for example, bonds, stocks and options). Risky assets evolution is described by Itô's equation. We assume also that there is a bank account where money grows exponentially according to the short rate (we shall call it risk-free rate). Evolution of the risky assets is described by Itô's equations; in particular, this setting covers a case when there are $m$ driving Brownian motions and $N >> m$ bonds which with different maturing times $T_1, \ldots, T_N$. The problem is to find an investment strategy for which $\mathbf{E}U(\widetilde{X}(T))$ is to be maximized, where $\mathbf{E}$ denotes the mathematical expectation, $U(\cdot)$ is an utility function, $\widetilde{X}(T) = \exp\left(-\int_0^T r(s)ds\right) X(T)$ is the normalized wealth, and where $X(T)$ represents the wealth at the final time $T$. There are many works devoted to different modifications of this problem (see, e.g., Merton (1969) and survey in Hakansson (1997) and Karatzas and Shreve (1998))). In the setting generally assumed in finance, cf. Merton (1990), Sec. 15.5, the coefficients are assumed to satisfy an Itô equation. Then the solution of the optimal investment problem can be obtained via dynamic programming approach. However, it is not easy to find the explicit solution by this method, because the corresponding Bellman equation is usually degenerate. Explicit formulas for optimal strategies have been obtained only for a few cases where appreciation rates are assumed to be non-random and known, and $U(\cdot)$ has quadratic form, log form or power form.

Investment problems for market where there are both bonds and stocks available are more difficult to study (some reasons for this were listed in Bielecki and Pliska (2001); in addition, we can add that the volatility matrix is not invertible in this case, and many standard methods are not applicable; moreover, the model descriptions there are usually cumbersome). However, the investment problems there were studied in dynamic programming approach in some cases, for example, with several driving Brownian motions (see *e.g.* Rutkowski (1997), Bielecki and Pliska (2001)).

We study the optimal investment problem for a diffusion market model such that the parameters $r(t)$, $a(t)$ and $\sigma(t)$ are all random; they are not adapted to the driving Brownian motion, but they are supposed to be currently observable (i.e. it is a case of "totally unhedgeable" coefficients, according to Karatzas and Shreve (1998), Chapter 6). In addition, we do not assume to know the distributions of $(r(\cdot), a(\cdot), \sigma(\cdot))$. Following Cvitanić and Karatzas (1999), and Cvitanić (2000), we consider the problem as a *maximin* problem: Find a strategy which maximizes the infimum of $\mathbf{E}U(\widetilde{X}(T))$ over all admissible $(r(\cdot), a(\cdot), \sigma(\cdot))$ from a given class; the process $(r(\cdot), a(\cdot), \sigma(\cdot))$ is supposed to be currently observable. For this problem, we show that the duality theorem holds under some non-restrictive conditions. Thus, the *maximin* problem which, as far as we know, cannot be solved directly, is effectively reduced to a *minimax* problem. Moreover, it is proved that *minimax* problem requires minimization only *over a single scalar parameter* $R$ even for multi-stock market, where $R = \int_0^T |\sigma(t)^{-1}(a(t) - r(t)\mathbf{1})|^2 dt$. This interesting effect follows from the result of Dokuchaev and Haussmann (2001) for the optimal compression problem. Using this effect, the original *maximin* problem is effectively solved; the optimal strategy is derived via solution of a *linear* parabolic equation.

Cvitanić and Karatzas (1999) and Cvitanić (2000) consider a related *minimax* and *maximin* problems of minimizing $\mathbf{E}(\xi_1 - X(T))^+$ subject to $X(T) \geq \xi_2$, where $\xi_1$ and $\xi_2$ are given claims, for similar admissible strategies which allow direct observations of appreciation rates (adapted to the driving Brownian motion); however, the maximization over parameters in the dual *minimax* problem was not reduced to the scalar minimization, and the solution was not given for the general case. Furthermore, we consider more general utility functions. Dokuchaev and Teo (1998) obtained a duality theorem for a problem in maximin setting with admissible strategies which use only historical prices.

## 2 Definitions and problem statement

We consider a market which consists of a risk free bond or bank account with price $B(t)$, $t \geq 0$, and $n$ risky stocks with prices $S_i(t)$, $t \geq 0$, $i = 1, 2, \ldots, n$, where $n < +\infty$ is given. The prices of the stocks evolve according to

$$dS_i(t) = S_i(t)\bigg(a_i(t)dt + \sum_{j=1}^m \sigma_{ij}(t)dw_j(t)\bigg), \quad t > 0, \tag{2.1}$$

where the $w_i(t)$ are standard independent Wiener processes, $a_i(t)$ are appreciation rates, and $\sigma_{ij}(t)$ are volatility coefficients. The initial price $S_i(0) > 0$ is a given nonrandom constant. The price of the bond evolves according to

$$B(t) = B(0) \exp\left(\int_0^t r(s)ds\right), \tag{2.2}$$

where $B(0)$ is a given constant that we take to be 1 without loss of generality, and $r(t)$ is the random process of the risk-free interest rate.

We assume that $w(\cdot) = (w_1(\cdot), \ldots, w_m(\cdot))$ is a standard Wiener process on a given standard probability space $(\Omega, \mathcal{F}, \mathbf{P})$, where $\Omega = \{\omega\}$ is a set of elementary events, $\mathcal{F}$ is a complete $\sigma$-algebra of events, and $\mathbf{P}$ is a probability measure.

We assume that $r(t)$, $a(t) \triangleq \{a_i(t)\}_{i=1}^n$, and $\sigma(t) \triangleq \{\sigma_{ij}(t)\}_{i,j=1}^{n,m}$ are uniformly bounded, measurable random processes.

We are interested in the case of *degenerate* $\sigma(t)\sigma(t)^\top$, because we want to cover, in particular, a case when the market includes $N_1$ stocks and $N = n - N_1$ zero-coupon bonds with different maturing times $T_1, \ldots, T_N$, where $N >> m$. We also want to cover the case when options on stocks and bonds are considered as tradable assets. Assumption 4.1 below will be in force throughout this paper and it ensures that the market is arbitrage free and at the same time allows us to include bonds and options into consideration (see, e.g., Lamberton and Lapeyre (1996)).

Set $\mu(t) \triangleq (r(t), \widetilde{a}(t), \sigma(t))$, where $\widetilde{a}(t) \triangleq a(t) - r(t)\mathbf{1}$.

Let $\{\mathcal{F}_t^\mu\}_{0 \leq t \leq T}$, be the filtration generated by the process $(S(t), \mu(t))$ completed with the null sets of $\mathcal{F}$. Clearly, $\mathcal{F}_t^\mu$ coincides with the filtration generated by the processes $(w(t), \mu(t))$, and with the filtration generated by the processes $(\widetilde{S}(t), \mu(t))$, where

$$p(t) \triangleq \exp\left(-\int_0^t r(s)ds\right) = B(t)^{-1}, \qquad \widetilde{S}(t) \triangleq p(t)S(t).$$

We describe now distributions of $\mu(\cdot)$ and what we suppose known about them.

We assume that there exist a finite-dimensional Euclidean space $\bar{E}$, a compact subset $\mathcal{T} \subset \bar{E}$, and a measurable function $M(\cdot) = (M_r(\cdot), M_a(\cdot), M_\sigma(\cdot)) : [0,T] \times \mathcal{T} \times C([0,T]; \mathbf{R}^n) \to \mathbf{R} \times \mathbf{R}^n \times \mathbf{R}^{n \times n}$ that is uniformly bounded and such that $M(t, \alpha, \xi)$ is continuous in $\alpha \in \mathcal{T}$ for all $t$ and $\xi \in C([0,t]; \mathbf{R})$ and

$$M(t, \cdot) = (M_r(t, \cdot), M_a(t, \cdot), M_\sigma(t, \cdot)) : \mathcal{T} \times C([0,t]; \mathbf{R}^n) \to \mathbf{R} \times \mathbf{R}^n \times \mathbf{R}^{n \times n},$$

We assume that $\mathcal{T}$ and $M(\cdot)$ are such that the solution of (2.1) with $\mu(t) = (r(t), \widetilde{a}(t), \sigma(t)) = M(t, \alpha, S(\cdot)|_{[0,t]})$ is well defined for any $\alpha \in \mathcal{T}$ as the unique strong solution of Itô's equation. Let $S_\alpha(\cdot)$ denote the corresponding solution.

For $\alpha \in \mathcal{T}$, set

$$\bar{M}_r(t, \alpha) \triangleq M_r\left(t, \alpha, S_\alpha(\cdot)|_{[0,t]}\right),$$
$$\bar{M}_a(t, \alpha) \triangleq M_a\left(t, \alpha, S_\alpha(\cdot)|_{[0,t]}\right),$$
$$\bar{M}_\sigma(t, \alpha) \triangleq M_\sigma\left(t, \alpha, S_\alpha(\cdot)|_{[0,t]}\right).$$

**Definition 2.1** *Let $\mathcal{A}(\mathcal{T})$ be a set of all random processes $\mu'(t) = (r'(t), \widetilde{a}'(t), \sigma'(t))$ such* *that there exists a random vector* $\Theta : \Omega \to \mathcal{T}$ *independent of* $w(\cdot)$ *and such that*

$$\begin{cases} r'(t) \equiv \bar{M}_r(t, \Theta) \\ \tilde{a}'(t) \equiv \bar{M}_a(t, \Theta) \\ \sigma'(t) \equiv \bar{M}_\sigma(t, \Theta). \end{cases} \quad (2.3)$$

*We assume that* $\mu(\cdot) \in \mathcal{A}(\mathcal{T})$, *and that is the only information available.* Notice that the solution of (2.1) is well defined for any $\mu(\cdot) \in \mathcal{A}(\mathcal{T})$, but the market is incomplete.

*Remark.* In fact, the solution of investment problem obtained below does not require to know $\bar{E}, \mathcal{T}$ and $M(\cdot)$.

Note that $\mathcal{T}$ can be interpreted as a set of unknown parameters.

**Example 2.1** Let $n = 1$, $\bar{E} = \mathbf{R}^N$, where $N > 0$ is an integer, $\mathcal{T} \subset E$ be a subset, and

$$(\bar{M}_r(t, \alpha), \bar{M}_\sigma(t, \alpha)) \equiv (r, \sigma), \qquad \alpha = (\alpha_1, \ldots, \alpha_N) \in \mathcal{T},$$
$$\bar{M}_a(t, \alpha) = \alpha_k, \quad t \in \left[\tfrac{(k-1)T}{N}, \tfrac{kT}{N}\right), \qquad k = 1, \ldots, N.$$

where $r$, $\sigma$ are constants. Then $\mathcal{A}(\mathcal{T})$ is the set of all processes $\mu(t) = (r(t), \tilde{a}(t), \sigma(t))$ such that

$$r(t) \equiv r, \quad \sigma(t) \equiv \sigma,$$
$$\tilde{a}(t) = \Theta_k, \quad t \in \left[\tfrac{(k-1)T}{N}, \tfrac{kT}{N}\right), \quad k = 1, \ldots, N,$$

where $\Theta = (\Theta_1, \ldots, \Theta_N)$ is a $N$-dimensional random vector independent of $w(\cdot)$, $\Theta \in \mathcal{T}$.

**Remark 2.1** It is easy to see that our description of the class of admissible $\mu(\cdot)$ covers a setting when the minimum of $R_\mu$ over the class is given, or when the class of admissible $\mu(\cdot)$ is defined by a condition $R_\mu \in [R_1, R_2]$, where $R_1, R_2$ are given, $0 \le R_1 < R_2 \le +\infty$. (It suffices to choose an appropriate pair $(\Theta, M(\cdot))$.)

For $\alpha \in \mathcal{T}$, set

$$\mu_\alpha(t) \triangleq (\bar{M}_r(t, \alpha), \bar{M}_a(t, \alpha), \bar{M}_\sigma(t, \alpha)),$$

where $\bar{M}_r(t, \alpha)$, $\bar{M}_a(t, \alpha)$ and $\bar{M}_\sigma(t, \alpha)$ are as in Definition 2.1.

Let $X_0 > 0$ be the initial wealth at time $t = 0$, and let $X(t)$ be the wealth at time $t > 0$, $X(0) = X_0$. We assume that

$$X(t) = \pi_0(t) + \sum_{i=1}^{n} \pi_i(t), \quad (2.4)$$

where the pair $(\pi_0(t), \pi(t))$ describes the portfolio at time $t$. The process $\pi_0(t)$ is the investment in the bond, $\pi_i(t)$ is the investment in the $i$th stock, $\pi(t) = (\pi_1(t), \ldots, \pi_n(t))^\top$, $t \geq 0$.

Let $\mathbf{S}(t) \triangleq \operatorname{diag}(S_1(t), \ldots, S_n(t))$ and $\widetilde{\mathbf{S}}(t) \triangleq \operatorname{diag}(\widetilde{S}_1(t), \ldots, \widetilde{S}_n(t))$ be diagonal matrices with the corresponding diagonal elements. The portfolio is said to be self-financing, if

$$dX(t) = \pi(t)^\top \mathbf{S}(t)^{-1} dS(t) + \pi_0(t) B(t)^{-1} dB(t). \tag{2.5}$$

It follows that for such portfolios

$$dX(t) = r(t) X(t) \, dt + \pi(t)^\top \left(\widetilde{a}(t) \, dt + \sigma(t) \, dw(t)\right), \tag{2.6}$$

$$\pi_0(t) = X(t) - \sum_{i=1}^n \pi_i(t),$$

so $\pi$ alone suffices to specify the portfolio; it is called a self-financing strategy.

The process $\widetilde{X}(t) \triangleq p(t) X(t)$ is called the normalized wealth. It satisfies

$$\widetilde{X}(t) = X(0) + \int_0^t p(t) \pi(s)^\top \widetilde{\mathbf{S}}(t)^{-1} d\widetilde{S}(s). \tag{2.7}$$

**Definition 2.2** *Let $\mathcal{G}_t$ be a filtration. Let $\widetilde{\Sigma}(\mathcal{G}.)$ be the class of all $\mathcal{G}_t$-adapted processes $\pi(\cdot) = (\pi_1(t), \ldots, \pi_n(t))$ such that*

$$\mathbf{E} \int_0^T |\pi(t)|^2 \, dt < \infty.$$

A process $\pi(\cdot) \in \widetilde{\Sigma}(\mathcal{G}.)$ is said to be an *admissible* strategy with corresponding wealth $X(\cdot)$.

Let $X_0 > 0$ be the initial wealth at time $t = 0$, and let $X(t)$ be the wealth at time $t > 0$. Let $\widetilde{X}(t)$ be the normalized wealth.

For an Euclidean space $E$ we shall denote by $B([0, T]; E)$ the set of bounded measurable functions $f(t) : [0, T] \to E$. By the definitions of $\widetilde{\Sigma}(\mathcal{F}^\mu)$ and $\mathcal{F}_t^\mu$, any admissible self-financing strategy is of the form

$$\pi(t) = \Gamma(t, [S(\cdot), \mu(\cdot)]|_{[0,t]}), \tag{2.8}$$

where $G(\cdot)$ is a measurable function, $\Gamma(t, \cdot) : B([0, t]; \mathbf{R}^n \times \mathbf{R} \times \mathbf{R}^n \times \mathbf{R}^{n \times n}) \to \mathbf{R}^n$, $t \geq 0$.

Clearly, the random processes $\pi(\cdot)$ with the same $\Gamma(\cdot)$ in (2.8) may be different for different $\mu(\cdot) = (r(\cdot), \widetilde{a}(\cdot), \sigma(\cdot))$. Hence we also introduce strategies defined by $\Gamma(\cdot)$: the function $\Gamma(\cdot)$ in (2.8) is said to be a CL-strategy (closed-loop strategy).

**Definition 2.3** *Let $\mathcal{C}$ be the class of all functions $\Gamma(t, \cdot) : B([0, t]; \mathbf{R}^n \times \mathbf{R} \times \mathbf{R}^n \times \mathbf{R}^{n \times n}) \to \mathbf{R}^n$, $t \geq 0$ such that the corresponding strategy $\pi(\cdot)$ defined by (2.8) belongs to $\widetilde{\Sigma}(\mathcal{F}_t^\mu)$ for any $\mu(\cdot) = (r(\cdot), \widetilde{a}(\cdot), \sigma(\cdot)) \in \mathcal{A}(\mathcal{T})$ and*

$$\sup_{\mu(\cdot) = \mu_\alpha(\cdot):\ \alpha \in \mathcal{T}} \mathbf{E} \int_0^T |\pi(t)|^2\, dt < \infty.$$

*A function $\Gamma(\cdot) \in \mathcal{C}$ is said to be an admissible CL-strategy.*

Let the initial wealth $X(0)$ be fixed. For an admissible self-financing strategy $\pi(\cdot)$ such that $\pi(t) = \Gamma(t, [S(\cdot), \mu(\cdot)]|_{[0,t]})$, the process $(\pi(t), X(t))$ is uniquely defined by $\Gamma(\cdot)$ and $\mu(\cdot) = (r(\cdot), \widetilde{a}(\cdot), \sigma(\cdot))$ given $w(\cdot)$. We shall use the notation $X(t, \Gamma(\cdot), \mu(\cdot))$ and $\widetilde{X}(t, \Gamma(\cdot), \mu(\cdot))$ to denote the corresponding total wealth and normalized wealth. Furthermore, we shall use the notation $S(t) = S(t, \mu(\cdot))$ and $\widetilde{S}(t) = \widetilde{S}(t, \mu(\cdot))$ to emphasize that the stock price is different for different $\mu(\cdot)$.

## 3 Problem statement

Let $T > 0$ and $X_0$ be given. Let $U(\cdot) : \mathbf{R} \to \mathbf{R} \cup \{-\infty\}$ be a given measurable function such that $U(X_0) < +\infty$. Let $\mathbf{D} \subset \mathbf{R}$ be a given convex set, $X_0 \in \mathbf{D}$.

We may state our general problem as follows: Find an admissible CL-strategy $\Gamma(\cdot)$ and the corresponding self-financing strategy $\pi(\cdot) \in \widetilde{\Sigma}(\mathcal{F}_t^\mu)$ that solves the following optimization problem:

$$\text{Maximize} \quad \min_{\mu(\cdot) \in \mathcal{A}(\mathcal{T})} \mathbf{E} U(\widetilde{X}(T, \Gamma(\cdot), \mu(\cdot))) \quad \text{over} \quad \Gamma(\cdot) \tag{3.1}$$

$$\text{subject to} \quad \begin{cases} X(0, \Gamma(\cdot), \mu(\cdot)) = X_0, \\ \widetilde{X}(T, \Gamma(\cdot), \mu(\cdot)) \in \mathbf{D} \quad \text{a.s.} \quad \forall \mu(\cdot) \in \mathcal{A}(\mathcal{T}). \end{cases} \tag{3.2}$$

Clearly, the maximin setting has no sense if, for example, $\mu(t) \equiv \Theta$, where $\Theta$ is a random element of $\mathbf{R} \times \mathbf{R}^n \times \mathbf{R}^{n \times n}$ which is constant in time; one can identify $\Theta$ instantly. However, the optimal solution for a more general case needs knowledge about distribution of future values of $\mu(\cdot)$.

**Example 3.1** *Let $n = 1$, $\mathcal{T} \triangleq \{\alpha_1, \alpha_2\}$, where $\alpha_i \in \mathbf{R}$. Let*

$$(\bar{M}_r(t, \alpha), \bar{M}_a(t, \alpha)) \equiv (r, \widetilde{a}), \qquad \bar{M}_\sigma(t, \alpha) = \begin{cases} \alpha_1, & t < T/2 \\ \alpha, & t \geq T/2, \end{cases}$$

where $r$, $\widetilde{a}$ are constants, i.e.,

$$(r(t), \widetilde{a}(t)) \equiv (r, \widetilde{a}), \qquad \sigma(t) = \begin{cases} \alpha_1, & t < T/2 \\ \Theta, & t \geq T/2, \end{cases}$$

where $\Theta$ is a random variable independent of $w(\cdot)$ which can have only two values, $\alpha_1$ and $\alpha_2$. Let $\kappa \in [0, 1)$ and $\mu(\cdot) \in \mathcal{A}(\mathcal{T})$ be given. Consider the problem

$$\text{Maximize} \quad \mathbf{E} \log \widetilde{X}(T, \Gamma(\cdot), \mu(\cdot)) \quad \text{over} \quad \Gamma(\cdot)$$

$$\text{subject to} \quad \begin{cases} X(0, \Gamma(\cdot), \mu(\cdot)) = X_0, \\ \widetilde{X}(T, \Gamma(\cdot), \mu(\cdot)) \geq \kappa X_0 \quad \text{a.s.} \end{cases}$$

By Theorem 5.1 (ii) from Dokuchaev and Haussmann (2001), it follows that if $\Theta \equiv \alpha_1$ or $\Theta \equiv \alpha_2$, then the optimal strategy exists, and if $\kappa \neq 0$, then the corresponding optimal strategies for these two cases differs at the time interval $[0, T/2)$ (see, e.g., Lemma 5.2 below). Hence the optimal strategy cannot be obtained from observations of historical $\widetilde{a}(t)$ and $S(t)$ without knowledge of future distributions. The only exception is the case $\kappa = 0$, when the optimal strategy given $\mu(\cdot)$ is myopic.

**The case of myopic strategies**

**Proposition 3.1** *Let $X_0 = X(0) > 0$ and let $\sigma(t)\sigma(t)^\top \geq cI_n$, where $c > 0$ is a constants and $I_n$ is the $n \times n$ identity matrix. Further, let $M(t, \alpha, \xi) \equiv M(t, \alpha)$ does not depend on $\xi \in C([0, T]; \mathbf{R}^n)$, and let one of the following conditions be satisfied:*

(i) $U(x) = \log(x)$, $\mathbf{D} = [0, +\infty)$;

(ii) $U(x) = x^\delta$, $\mathbf{D} = [0, +\infty)$, *where* $\delta < 1$, $\delta \neq 0$;

(iii) $U(x) = -kx^2 + cx$, $\mathbf{D} = \mathbf{R}$, *where* $k \in \mathbf{R}$ *and* $c \geq 0$.

*Then there exists $C_0, C_1, \nu \in \mathbf{R}$ such that $C_1 \neq 0$, $\nu \neq 0$ are constants, and that the optimal strategy $\pi(\cdot) \in \widetilde{\Sigma}(\mathcal{F}^\mu_\cdot)$ for the problem (3.1)-(3.2) has the form*

$$\pi(t)^\top = \nu B(t)(\widetilde{X}(t) - C_0)\widetilde{a}(t)^\top Q(t), \tag{3.3}$$

*where $\widetilde{X}(t)$ is the corresponding normalized wealth, $Q(t) \triangleq (\sigma(t)\sigma(t)^\top)^{-1}$. This solution is optimal for the problem*

$$\text{Maximize} \quad \mathbf{E} U(\widetilde{X}(T, \Gamma(\cdot), \mu(\cdot))) \quad \text{over} \quad \Gamma(\cdot) \tag{3.4}$$

*for any $\mu(\cdot)$.*

# 4 Additional assumptions and some examples

We assume that $\mu(\cdot) \in \mathcal{A}(\mathcal{T})$.

To proceed further, we assume that the following Conditions 4.1-4.6 remain in force throughout this paper. The first of them ensures that the market is arbitrage free.

**Condition 4.1** *For any $\mu(\cdot) \in \mathcal{A}(\mathcal{T})$, there exists a random process $\theta_\mu(t) = (\theta_{\mu 1}(t), \ldots, \theta_{\mu m}(t))^\top$ such that*

$$\int_0^T |\theta_\mu(t)|^2 dt < +\infty, \quad \sigma(t)\theta_\mu(t) \equiv \widetilde{a}(t). \tag{4.1}$$

In addition, we assume without a loss of generality that there exists a set $\{i_1, \ldots, i_m\} \subseteq \{1, \ldots, n\}$ such that the matrix $\sigma'(t) \triangleq \{\sigma_{i_k j}(t)\}_{k,j=1}^m$ is such that $\sigma'(t)\sigma'(t)^\top \geq c' I_m$, a.s., $\forall t$, where $c' > 0$ is a constants and $I_m$ is the $m \times m$ identity matrix.

Note that the process $\theta_\mu(\cdot)$ is uniquely defined given $\mu(\cdot)$.

Set
$$R_\mu \triangleq \int_0^T |\theta_\mu(t)|^2 dt.$$

Set
$$\mathcal{Z}(t, \mu(\cdot)) \triangleq \exp\left(\int_0^t \theta_\mu(s)^\top dw(s) + \frac{1}{2}\int_0^t |\theta_\mu(s)|^2 ds\right). \tag{4.2}$$

Our standing assumptions imply that $\mathbf{E}\mathcal{Z}(T, \mu_\alpha(\cdot))^{-1} = 1$ for all $\alpha \in \mathcal{T}$.

Define the (equivalent martingale) probability measure $\mathbf{P}_*^\alpha$ by
$$\frac{d\mathbf{P}_*^\alpha}{d\mathbf{P}} = \mathcal{Z}(T, \mu_\alpha(\cdot))^{-1}.$$

Let $\mathbf{E}_*^\alpha$ be the corresponding expectation.

**Condition 4.2** *There exists a measurable set $\Lambda \subseteq \mathbf{R}$, and a measurable function $F(\cdot, \cdot) : (0, \infty) \times \Lambda \to \mathbf{D}$ such that for each $z > 0$, $\widehat{x} = F(z, \lambda)$ is a solution of the optimization problem*

$$\text{Maximize} \quad zU(x) - \lambda x \quad \text{over } x \in \mathbf{D}. \tag{4.3}$$

*Moreover, this solution is unique for a.e. $z > 0$.*

**Condition 4.3** *For any $\alpha \in \mathcal{T}$, there exist $\widehat{\lambda}_\alpha \in \Lambda$, $C = C_\alpha > 0$, and $c_0 = c_{0,\alpha} \in (0, 1/(2R_{\mu_\alpha}))$ such that $F(\cdot, \widehat{\lambda})$ is piecewise continuous on $(0, \infty)$, $F(\mathcal{Z}(T, \mu_\alpha(\cdot)), \widehat{\lambda}_\alpha)$ is $\mathbf{P}_*^\alpha$-integrable, and*

$$\begin{cases} \mathbf{E}_*^\alpha F(\mathcal{Z}(T, \mu_\alpha(\cdot)), \widehat{\lambda}_\alpha) = X_0, \\ |F(z, \widehat{\lambda}_\alpha)| \leq C z^{c_0 \log z} \quad \forall z > 0. \end{cases} \tag{4.4}$$

Some examples when conditions similar to the imposed above ones are satisfied can be found in Dokuchaev and Haussmann (2001).

**Condition 4.4** *The function $U(x) : \mathbf{R} \to \mathbf{R}$ is either concave or convex in $x \in \mathbf{D}$, and there exist constants $c > 0$, $p \in (1, 2]$, $q \in (0, 1]$ such that*

$$\begin{aligned}|U(x)| &\leq c\left(|x|^p + 1\right), \\ |U(x) - U(x_1)| &\leq c\left(1 + |x| + |x_1|\right)^{2-q} |x - x_1|^q \quad \forall x, x_1 \in \mathbf{D}.\end{aligned} \quad (4.5)$$

Notice that condition 4.5 is not restrictive if $\mathbf{D} \subset \mathbf{R}$ is a bounded interval (the case that is not excluded; this case includes goal achieving problem as well as any problem where an investor wish to avoid big variance for sure).

**Condition 4.5** (i) *The set $\mathcal{T}$ is such that $\mathcal{T} = \mathcal{T}^{(1)} \times \mathcal{T}^{(2)}$, where the set $\mathcal{T}^{(2)}$ is either finite or countable, i.e., $\mathcal{T} = \{\alpha\} = \{(\alpha^{(1)}, \alpha^{(2)})\}$ and $\mathcal{T}^{(2)} = \{\alpha_1^{(2)}, \alpha_2^{(2)}, \ldots\}$;*

(ii) *the function $(M_r(t, \alpha, \xi), M_\sigma(t, \alpha, \xi))$, where $\alpha = (\alpha^{(1)}, \alpha^{(2)})$, does not depend on $\alpha^{(1)}$ given $\alpha^{(2)}$ and $\xi \in C([0, t]; \mathbf{R}^n)$, i.e., $M_r(t, \alpha, \xi) \equiv M_r(t, \alpha^{(2)}, \xi)$ and $M_\sigma(t, \alpha, \xi) \equiv M_\sigma(t, \alpha^{(2)}, \xi)$; and*

(iii) *If the function $M_a(t, \alpha, \xi)$, where $\alpha = (\alpha^{(1)}, \alpha^{(2)})$, does depend on $\alpha^{(1)}$ given $\alpha^{(2)}$ then Condition 4.4 is satisfied with $p \in (1, 2)$.*

In other words, the diffusion may depend only on discrete random variable independent on Wiener process. Notice that Condition 4.5 looks restrictive, but in fact it is rather technical, since the total number of elements of $\mathcal{T}^{(2)}$ may be unbounded. In particular, this condition is always satisfied when $p < 2$ and $\sigma(t)$ does not depend on $\Theta$; or when $\mathcal{T}^{(2)} = \emptyset$, i.e. $\mathcal{T} = \mathcal{T}^{(2)}$.

**Condition 4.6** *At least one of the following conditions holds:*

(i) *$(\tilde{a}(t), \sigma(t)) \equiv (M_a(t, \Theta), M_\sigma(t, \Theta))$ with a deterministic function $(M_a(\cdot), M_\sigma(\cdot)) : [0, T] \times \mathcal{T} \to \mathbf{R}^n \times \mathbf{R}^{n \times n}$;*

(ii) *The matrix $\sigma(t)$ is diagonal for all arguments, and*

$$\sigma(t) \equiv M_\sigma(t, \Theta), \quad \tilde{a}_i(t) = \xi_i(t, \Theta, S(\cdot)|_{[0,t]}) \eta_i(t, \Theta),$$

*where $M_\sigma(\cdot) : [0, T] \times \mathcal{T} \to \mathbf{R}^{n \times n}$, $\xi_i(t, \cdot) : \mathcal{T} \times C([0, t]; \mathbf{R}^n) \to \mathbf{R}$ and $\eta_i(\cdot) : [0, T] \times \mathcal{T} \to \mathbf{R}$ are deterministic measurable functions such that $|\xi_i(t, \Theta, S(\cdot)|_{[0,t]})| \equiv 1$, $i = 1, \ldots n$.*

**Some examples**

**Example 4.1** *(Multi-bond market).* Consider a market with zero-coupon bonds with prices $P(t,T_k)$, where $t \leq T_k$, and where $\{T_k\}_{k=1}^N$ is a given set of maturing times, $0 < T_1 < \cdots < T_N = T$, $P(T_k, T_k) = 1$. We assume that investor can buy and sell bonds on this market. Let $\mathcal{F}_t^1$ be a filtration generated by the scalar Wiener process $w_1(t)$. Let $P(t, T_k)$ be adapted to $\mathcal{F}_t^1$, and let $B(t)$ be the "risk-free" asset such as defined above with risk-free rate $r(t)$ adapted to $\mathcal{F}_t^1$. It is shown in cf. Lamberton and Lapeyre (1996), section 6.2.1, that if this bond market is arbitrage free then there exists a $\mathcal{F}_t^1$-adapted process $q(t)$ such that

$$P(t, T_k) = \mathbf{E}\left\{\exp\left(-\int_t^{T_k} r(s)ds + \int_t^{T_k} q(s) dw_1(s) - \tfrac{1}{2}\int_t^{T_k} q(s)^2 ds\right) \middle| \mathcal{F}_t^1\right\}. \qquad (4.6)$$

On the other hand, under some mild conditions, any $\mathcal{F}_t^1$-adapted process $q(t)$ defines an arbitrage free bond market with prices (4.6).

By Proposition 6.1.3 from Lamberton and Lapeyre (1996), for any $T_k$, there exists a $\mathcal{F}_t^1$–adapted in $t$ process $\sigma(t, T_k)$ such that

$$d_t P(t, T_k) = P(t, t_k)\Big([r(t) - q(t)\sigma(t, T_k)]\, dt + \sigma(t, T_k) dw(t)\Big), \quad t < T_k.$$

Then we can treat this market as a special case of our market, where $m = 1$, $n = N$, and the set of risky assets is $S_k(t) = P(t, T_k)$, $k = 1, \ldots, N$, and $\mu(t) = (r(t), \widetilde{a}(t), \sigma(t))$, where $\widetilde{a}(t) \equiv (\widetilde{a}_1(t), \ldots, \widetilde{a}_k(t))^\top \in \mathbf{R}^n$, $\sigma(t) \equiv (\sigma_{11}(t), \ldots, \sigma_{k1}(t))^\top \in \mathbf{R}^{n \times 1}$, and where

$$(\widetilde{a}_k(t), \sigma_{k1}(t)) = \begin{cases} (-q(t)\sigma(t, T_k), \sigma(t, T_k)), & t \leq T_k \\ (0, 0), & t > T_k. \end{cases}$$

Condition 4.1 is satisfied with $\theta_\mu(t) \equiv -q(t)$.

Note that a special case when $q(t)$ is a deterministic process is a modification of the Vasicec model, where $q(t)$ is a constant (see Lamerton and Lapeyre (1996), p.127). If $q(t)$ is deterministic, then $R_\mu$ is non-random.

Clearly, this generic model can be easily developed further for a model that contains $m > 1$ driving Brownion motions, and contains both stocks and bonds.

**Example 4.2** *(Stock and options market)* Consider a risky asset (stock) $S_1(t)$ defined by (2.1) with $i = 1$ and $m = 1$. Let $r(t) \equiv r$ and $\sigma_{11}(t) \equiv \sigma_{11} \neq 0$ be given non-random constants. Further, we assume that there are available European options on that stocks with the same expiration time $T$ and different strike prices $K_1, \ldots, K_N$, where $N$ is an integer, possibly a large number. Let $H_{BS,c}(t, x, K)$ and $H_{BS,p}(t, x, K)$ denote Black-Scholes prices for the put and call options with the claims $(S_1(T) - K)^+$ and $(K - S_1(T))^+$ respectively given condition $\widetilde{S}_1(t) = x$, where $\widetilde{S}_1(t) \triangleq e^{-rt} S_1(t)$. We shall consider options as additional tradable assets, i.e. we shall consider *stock-options market*. Then we can treat this market as a special case of our market, where $m = 1$, $n = 1 + 2N$, and where the normalized prices for the risky assets are $\widetilde{S}_1(t), \ldots, \widetilde{S}_n(t)$, where

$$\widetilde{S}_{1+i}(t) = \widetilde{H}_{BS,c}(t, \widetilde{S}(t), K_i), \quad \widetilde{S}_{1+N+i}(t) \triangleq \widetilde{H}_{BS,p}(t, \widetilde{S}(t), K_i), \quad i = 1, \ldots, N$$
$$\widetilde{H}_{BS,c}(t, x, K) \triangleq e^{-rt} H_{BS,c}(t, x, K), \quad \widetilde{H}_{BS,c}(t, x, K) \triangleq e^{-rt} H_{BS,c}(t, x, K).$$

The well known Black-Scholes formula for the option prices prices ensure that

$$d\widetilde{S}_{1+i}(t) = \frac{\partial \widetilde{H}_{BS,c}}{\partial x}(t, \widetilde{S}_1(t), K_i) d\widetilde{S}_1(t),$$
$$d\widetilde{S}_{1+N+i}(t) = \frac{\partial \widetilde{H}_{BS,p}}{\partial x}(t, \widetilde{S}_1(t), K_i) d\widetilde{S}_1(t), \quad i = 1, \ldots, N.$$

Then $\mu(t) = (r, \widetilde{a}(t), \sigma(t))$, where $\widetilde{a}(t) = (\widetilde{a}_1(t), \ldots, \widetilde{a}_n(t))$, $\sigma(t) \equiv (\sigma_{11}(t), \ldots, \sigma_{k1}(t))^\top \in \mathbf{R}^{n \times 1}$, and where

$$(\widetilde{a}_k(t), \sigma_{k1}(t)) = \begin{cases} (\widetilde{a}_1(t), \sigma_{11}), & k = 1 \\ \left( \widetilde{a}_1(t) \frac{\partial \widetilde{H}_{BS,c}}{\partial x}(t, \widetilde{S}_1(t), K_i), \sigma_{11} \frac{\partial \widetilde{H}_{BS,c}}{\partial x}(t, \widetilde{S}_1(t), K_i) \right), & 1 < k \leq 1 + N \\ \left( \widetilde{a}_1(t) \frac{\partial \widetilde{H}_{BS,c}}{\partial x}(t, \widetilde{S}_1(t), K_i), \sigma_{11} \frac{\partial H_{BS,p}}{\partial x}(t, \widetilde{S}_1(t), K_i) \right), & k > 1 + N. \end{cases}$$

Assumption 4.1 is satisfied with $\theta_\mu(t) \equiv \sigma_{11}^{-1} \widetilde{a}_1(t)$.

**Example 4.3** *(Random time of volatility change)* Let $k > 0$ be a integer. Let $n = 1$,

$$(r(t), \widetilde{a}(t)) \equiv (r, \widetilde{a}), \qquad \sigma(t) = \begin{cases} \bar{\sigma}, & t < \tau \\ \Theta, & t \geq \tau, \end{cases}$$

where $r > 0$, $\widetilde{a}$, and $\bar{\sigma}$ are constants, $\tau$ and $\Theta$ are random variables such that the pair $(\tau, \Theta)$ is independent of $w(\cdot)$, and such that

$$\Theta \in \left\{ 0, \pm \tfrac{1}{k}, \ldots, \pm \tfrac{k-1}{k}, \pm 1 \right\}, \qquad \tau \in \left\{ 0, \tfrac{1}{k}, \ldots, \tfrac{k-1}{k}, 1 \right\}.$$

Then Condition 4.5(i) is satisfied with $\bar{E} = \mathbf{R}^2$,

$$\mathcal{T} = \left\{ \alpha = (\alpha_1, \alpha_2) : \alpha_1 \in \left\{ 0, \pm \tfrac{1}{k}, \ldots, \pm \tfrac{k-1}{k}, 1 \right\}, \alpha_2 \in \left\{ 0, \tfrac{1}{k}, \ldots, \tfrac{k-1}{k}, 1 \right\} \right\},$$
$$(M_r(t, \alpha), M_a(t, \alpha)) \equiv (r, \widetilde{a}), \qquad M_\sigma(t, \alpha) = \begin{cases} \bar{\sigma}, & t < \alpha_2 \\ \alpha_1, & t \geq \alpha_2. \end{cases}$$

## 5 The main result: solution of the maximin problem

For given $R > 0$, $\lambda \in \Lambda$, let the function $H(\cdot) = H(\cdot, R, \lambda) : \mathbf{R}_+ \times [0, T] \to \mathbf{R}$ be the solution of the following Cauchy problem:

$$\begin{cases} \frac{\partial H}{\partial t}(x, t, R, \lambda) + \frac{R}{2T} x^2 \frac{\partial^2 H}{\partial x^2}(x, t, R, \lambda) = 0, \\ H(x, T, R, \lambda) = F(x, \lambda), \end{cases} \quad (5.1)$$

where $F(\cdot)$ is defined in Condition 4.2.

Introduce a function $\widetilde{\Gamma}(t, \cdot) : B([0, t]; \mathbf{R}^n \times \mathbf{R} \times \mathbf{R}^n \times \mathbf{R}^{n \times n}) \times (0, +\infty) \times \Lambda \to \mathbf{R}^n$ such that

$$\widetilde{\Gamma}(t, [S(\cdot), \mu(\cdot)]|_{[0,t]}, R, \lambda) = B(t) \frac{\partial H}{\partial x} \left[ \mathcal{Z}(t, \mu(\cdot)), \tau_\mu(t, R), R, \lambda \right] \mathcal{Z}(t, \mu(\cdot)) \widetilde{a}(t)^\top Q(t),$$

where the process $\mathcal{Z}(t, \mu(\cdot))$ is defined by (4.2) and where

$$\tau_\mu(t, R) = \tau(t, [S(\cdot), \mu(\cdot)]|_{[0,t]}, R) \triangleq \frac{T}{R} \int_0^t |\theta_\mu(s)|^2 ds.$$

Further, for a given $\alpha \in \mathcal{T}$, $R \geq 0$, let CL-strategy $\widehat{\Gamma}_\alpha(\cdot, R)$ be defined as

$$\widehat{\Gamma}_\alpha(t, [S(\cdot), \mu(\cdot)]|_{[0,t]}, R) \triangleq \begin{cases} \widetilde{\Gamma}(t, [S(\cdot), \mu(\cdot)]|_{[0,t]}, R, \widehat{\lambda}_\alpha) & \text{if } R > 0 \\ 0 & \text{if } R = 0, \end{cases}$$

where $\widehat{\lambda}_a$ is defined from Condition 4.3.

Note that Condition 4.6 ensures that $\mathbf{R}_{\mu_\alpha}$ is deterministic for any $\alpha \in \mathcal{T}$ is satisfied.

**Definition 5.1** *Let $\mathcal{C}_0$ be the set of all admissible CL-strategies $\Gamma(\cdot) \in \mathcal{C}$ such that*

$$\widetilde{X}(T, \Gamma(\cdot), \mu(\cdot)) \in \mathbf{D} \quad a.s. \quad \forall \mu(\cdot) \in \mathcal{A}(\mathcal{T}).$$

**Lemma 5.1** *For any $\mu(\cdot) = (r(\cdot), \widetilde{a}(\cdot), \sigma(\cdot))$, there exists $n \times m$-dimensional $\mathcal{F}_t$-adapted matrix process $D_\mu(t)$ such that*

$$\begin{cases} \theta_\mu(t)^\top D_\mu(t) \sigma(t) \equiv \theta_\mu(t)^\top \\ D_\mu(t) \widetilde{a}(t) \equiv \theta_\mu(t). \end{cases} \quad (5.2)$$

To formulate our main result, we shall need some generalizations of results from Dokuchaev and Haussmann (2001) for our market when the matrix $\sigma(t) \sigma(t)^\top$ can be degenerate, and these are summarized in the following lemma.

**Lemma 5.2** *(i) For any $R > 0$, $\lambda \in \Lambda$, the problem (5.1) has a unique solution $H(\cdot, R, \lambda) \in C^{2,1}((0, \infty) \times (0, T))$, with $H(x, t, R, \lambda) \to F(x, \lambda)$ a.e. as $t \to T-$.*

*(ii) For any $\alpha \in \mathcal{T}$, the strategy*

$$\widehat{\Gamma}_\alpha(t, [S(\cdot), \mu(\cdot)]|_{[0,t]}, R_{\mu_\alpha}) = B(t)\frac{\partial H}{\partial x}\left[\mathcal{Z}(t, \mu(\cdot)), \tau_\mu(t, R_{\mu_\alpha}), R_{\mu_\alpha}, \widehat{\lambda}_\alpha\right] \mathcal{Z}(t, \mu(\cdot))\widetilde{\theta}_\mu(t)^\top D_\mu(t) \tag{5.3}$$

*belongs to $\mathcal{C}_0$ and*

$$\mathbf{E}U(\widetilde{X}(T, \widehat{\Gamma}_\alpha(\cdot, R_{\mu_\alpha}), \mu_\alpha(\cdot))) \geq \mathbf{E}U(\widetilde{X}(T, \Gamma(\cdot), \mu_\alpha(\cdot))) \quad \forall \Gamma(\cdot) \in \mathcal{C}_0, \ \forall \alpha \in \mathcal{T}. \tag{5.4}$$

*(iii) The functions $F(\cdot, \widehat{\lambda}_\alpha)$, $H(\cdot, R_{\mu_\alpha}, \widehat{\lambda}_\alpha)$, $\widehat{\Gamma}_\alpha(\cdot, R_{\mu_\alpha})$ as well as the probability distribution of the optimal normalized wealth $\widetilde{X}(T, \widehat{\Gamma}_\alpha(\cdot), \mu_\alpha(\cdot))$ is uniquely defined by $R_{\mu_\alpha}$.*

*(iv) Let $\alpha_i \in \mathcal{T}$, $i = 1, 2$ be such that $R_{\mu_1} < R_{\mu_2}$, where $\mu_i \triangleq \mu_{\alpha_i}$. Then*

$$\mathbf{E}U(\widetilde{X}(T, \widehat{\Gamma}_{\alpha_1}(\cdot, R_{\mu_1}), \mu_1(\cdot))) < \mathbf{E}U(\widetilde{X}(T, \widehat{\Gamma}_{\alpha_2}(\cdot, R_{\mu_2}), \mu_2(\cdot))).$$

Set

$$R_{min} \triangleq \inf_{\mu(\cdot) \in \mathcal{A}(\mathcal{T})} R_\mu. \tag{5.5}$$

By the assumptions, $R_{min}$ is supposed to be known.

**Theorem 5.1** *(i) If $R_{min} = 0$, then the trivial strategy, $\Gamma(\cdot) \equiv 0$, is the unique optimal strategy in $\mathcal{C}$ for the problem (3.1)-(3.2).*

*(ii) Let $R_{min} > 0$, and let $\widehat{\alpha} \in \mathcal{T}$ be such that $R_{\widehat{\mu}} = R_{min}$, where $\widehat{\mu} \triangleq \mu_{\widehat{\alpha}}$. Then the strategy*

$$\widehat{\Gamma}_{\widehat{\alpha}}(t, [S(\cdot), \mu(\cdot)]|_{[0,t]}, R_{\widehat{\mu}}) \triangleq \widetilde{\Gamma}(t, [S(\cdot), \mu(\cdot)]|_{[0,t]}, R_{min}, \widehat{\lambda}_{\widehat{\alpha}}) \tag{5.6}$$

*belongs to $\mathcal{C}_0$ and is optimal in $\mathcal{C}$ for the problem (3.1)-(3.2).*

**Corollary 5.1** *The optimal strategy for the problem (3.1)-(3.2) does not depend on $(\mathcal{T}, \mathcal{M}(\cdot))$, if $R_{min} > 0$ is fixed and known.*

# 6  Appendix: Proofs

*Proof of Proposition 3.1.* By the assumptions on $M(\cdot)$, it follows that $\mu(\cdot) \in \mathcal{A}(\mathcal{T})$ does not depend on $w(\cdot)$. Condition 4.2 is satisfied with $F(x, \lambda) = C_1 \left(\frac{x}{\lambda}\right)^\nu + C_0$, where $C_1 \neq 0$, $C_0$ and $\nu \neq 0$ are constants. Then the proof follows from Corollary 5.1 from Dokuchaev and Haussmann (2001). □

*Proof of Lemma 5.1.* Without a loss of generality, we assume that the matrix $\widehat{\sigma}(t) \triangleq \{\sigma_{i,j}(t)\}_{i,j=1}^{m}$ is such that $\widehat{\sigma}(t)\widehat{\sigma}(t)^{\top} \geq cI_m$, a.s., $\forall t$, where $c > 0$ is a constants and $I_n$ is the $n \times n$ identity matrix. Set

$$D_\mu(t) \triangleq \left(\widehat{\sigma}(t)^{-1};\ 0_{n-m,m}\right),$$

where $0_{n-m,m}$ is the nil matrix in $\mathbf{R}^{n-m,m}$. Clearly, $D_\mu(t)\sigma(t) \equiv I_m$, then the first equation in (5.2) is satisfied. Further, let $m$-dimensional vector process $\widehat{a}(t)$ be such that $\widehat{a}_i(t) \triangleq \widetilde{a}_i(t)$, $i = 1, \ldots, m$. By the definition of $\theta$, we have that $\widehat{a}(t) = \widehat{\sigma}(t)\theta(t)$. Hence

$$D_\mu(t)\widetilde{a}(t) = \widehat{\sigma}(t)^{-1}\widehat{a}(t) = \widehat{\sigma}(t)^{-1}\widehat{\sigma}(t)\theta_\mu(t) = \theta_\mu(t).$$

This completes the proof. □

*Proof of Lemma 5.2.* Let $n = m$. Then statements (i)–(iii) follow immediately from Lemma 4.1, Theorem 5.1, and Lemma A.2 from Dokuchaev and Haussmann (2001). Let us show that statement (iv) holds. Let $\alpha_1 \in \mathcal{T}$ and $\alpha_2 \in \mathcal{T}$ be such that $R_{\mu_1} < R_{\mu_2}$, where $\mu_i \triangleq \mu_{\alpha_i}$. Further, let $\widehat{\mu}_{\alpha_2}(\cdot)$ be a process that is independent of $(\mu_{\alpha_1}(\cdot), w(\cdot))$ and has the same distribution as $\mu_{\alpha_2}(\cdot) \in \mathcal{A}(\mathcal{T})$. Consider a new auxiliary market with $2n$ stocks that consists of two independent groups of stocks that correspond to $\mu_{\alpha_1}(\cdot)$ and $\widehat{\mu}_{\alpha_2}(\cdot)$ (their driving Brownian motions and coefficients are independent). Then statement (iv) is a special case of Theorem 6.1 from Dokuchaev and Haussmann (2001), applied for the new market.

Let $n > m$. Then, similarly Dokuchaev and Haussmann (2001), it can be seen that $F(\mathcal{Z}(T, \mu_\alpha(\cdot), \widehat{\lambda}_\alpha)$ is the optimal claim, and this claim can be replicated by the strategy (5.3). This completes the proof. □

## A.1 Additional definitions

Without loss of generality, we describe the probability space as follows: $\Omega = \mathcal{T} \times \Omega'$, where $\Omega' = C([0, T]; \mathbf{R}^n)$. We are given a $\sigma$-algebra $\mathcal{F}'$ of subsets of $\Omega'$ generated by cylindrical sets, and a $\sigma$-additive probability measure $\mathbf{P}'$ on $\mathcal{F}'$ generated by $w(\cdot)$. Furthermore, let $\mathcal{F}_\mathcal{T}$ be the $\sigma$-algebra of all Borel subsets of $\mathcal{T}$, and $\mathcal{F} = \mathcal{F}_\mathcal{T} \otimes \mathcal{F}'$. We assume also that each $\mu(\cdot) \in \mathcal{A}(\mathcal{T})$ generates the $\sigma$-additive probability measure $\nu_\mu$ on $\mathcal{F}_\mathcal{T}$ (this measure is generated by $\Theta$ which corresponds to $\mu(\cdot)$).

Let $\widetilde{\Sigma}^R(\mathcal{F}^\mu_\cdot)$ be the enlargement of $\widetilde{\Sigma}(\mathcal{F}^\mu_\cdot)$ produced by replacing the filtration $\mathcal{F}^\mu_t$ by $\mathcal{F}^R_t$ generated by $\mathcal{F}^\mu_t$ and $R_\mu$ in the definition. (Note that the corresponding strategies are not adapted to $\mu(t)$.)

By the definitions of $\widetilde{\Sigma}^R(\mathcal{F}^\mu_\cdot)$, any admissible self-financing strategy from this class is of the form
$$\pi(t) = \Gamma(t, [S(\cdot), \mu(\cdot)]|_{[0,t]}, R_\mu),$$
where $\Gamma(t, \cdot) : B([0,t]; \mathbf{R} \times \mathbf{R}^n \times \mathbf{R}^N) \times \mathbf{R} \to \mathbf{R}^n$ is a measurable function, $t \geq 0$. With $\mathbf{P}(\cdot)$ replaced by $\mathbf{P}(\cdot \mid \mathcal{F}^R_0)$, we may apply Lemma 5.2 to obtain the optimal $\pi$ in the class $\widetilde{\Sigma}^R(\mathcal{F}^\mu_\cdot)$ for any $\mu(\cdot) \in \mathcal{A}(\mathcal{T})$. Note that the optimal strategy depends on random $R_\mu$, i.e. on $\Theta$; by Condition 4.6, there exists a measurable function $\phi : \mathcal{T} \to \mathbf{R}$ such that $R_\mu = \phi(\Theta)$.

Let $\mathring{\mathbf{R}}^n_+ \triangleq (0, +\infty)^n$.

For a function $\Gamma(t, \cdot) : C([0,t]; \mathring{\mathbf{R}}^n_+) \times B([0,t]; \mathbf{R} \times \mathbf{R}^n \times \mathbf{R}^{n \times n}) \times \mathbf{R}_+ \to \mathbf{R}^n$, introduce the following norm:
$$\|\Gamma(\cdot)\|_{\mathbf{X}} \triangleq \sup_{\mu(\cdot) = \mu_\alpha(\cdot):\ \alpha \in \mathcal{T}} \left( \sum_{i=1}^n \int_0^T \Gamma_i(t, [S(\cdot), \mu(\cdot)]|_{[0,t]}, R_\mu)^2 dt \right)^{1/2}. \tag{A.1}$$

**Definition A.1** *Let $\mathcal{C}^R_0$ be the set of all admissible CL-strategies $\Gamma(t, \cdot) : B([0,t]; \mathbf{R} \times \mathbf{R}^n \times \mathbf{R}^N) \times \mathbf{R} \to \mathbf{R}^n$ such that $\pi(t) = \Gamma(t, [S(\cdot), \mu(\cdot)]|_{[0,t]}, R_\mu) \in \widetilde{\Sigma}^R(\mathcal{F}^\mu_\cdot)$ for any $\mu(\cdot) \in \mathcal{A}(\mathcal{T})$, $\|\Gamma(\cdot)\|_{\mathbf{X}} < +\infty$ and*
$$\widetilde{X}(T, \Gamma(\cdot), \mu(\cdot)) \in \mathbf{D} \quad a.s. \quad \forall \mu(\cdot) \in \mathcal{A}(\mathcal{T}).$$

In fact, $\mathcal{C}^R_0$ is a subset of a linear space of functions with the norm (A.1).

## A.2 A duality theorem

To prove Theorem 5.1, we need the following duality theorem.

**Theorem A.1** *The following holds:*
$$\sup_{\Gamma(\cdot) \in \mathcal{C}^R_0} \inf_{\mu(\cdot) \in \mathcal{A}(\mathcal{T})} \mathbf{E} U(\widetilde{X}(T, \Gamma(\cdot), \mu(\cdot)))$$
$$= \inf_{\mu(\cdot) \in \mathcal{A}(\mathcal{T})} \sup_{\Gamma(\cdot) \in \mathcal{C}^R_0} \mathbf{E} U(\widetilde{X}(T, \Gamma(\cdot), \mu(\cdot))). \tag{A.2}$$

To prove Theorem A.1, we need several preliminary results, which are presented below as lemmas. The first of which is

**Lemma A.1** *The function $\widetilde{X}(T, \Gamma(\cdot), \mu(\cdot))$ is affine in $\Gamma(\cdot)$.*

*Proof.* By (2.7), it follows that $\widetilde{X}(t) = \widetilde{X}(t, \Gamma(\cdot), \mu(\cdot))$ satisfies

$$\widetilde{X}(t) = X(0) + \sum_{i=1}^{n} \int_0^t p(\tau) \Gamma_i(\tau, [S(\cdot, \mu(\cdot)), \mu(\cdot)]_{[0,\tau]}, R_\mu) \Big( \tilde{a}_i(t) dt + \sum_{j=1}^{n} \sigma_{ij}(t) dw_j(\tau) \Big). \tag{A.3}$$

It is easy to see that $\widetilde{X}(T, \Gamma(\cdot), \mu(\cdot))$ is affine in $\Gamma(\cdot)$. This completes the proof. □

**Lemma A.2** *The set $\mathcal{C}_0^R$ is convex.*

*Proof.* Let $p \in (0,1)$, $\mu(\cdot) \in \mathcal{A}(\mathcal{T})$, $\Gamma^{(i)}(\cdot) \in \mathcal{C}_0^R$, $i = 1, 2$, and

$$\Gamma(\cdot) \triangleq (1-p)\Gamma^{(i)}(t) + p\Gamma^{(i)}(\cdot).$$

By Lemma A.1, it follows that

$$\widetilde{X}(T, \Gamma(\cdot), \mu(\cdot)) = (1-p)\widetilde{X}(T, \Gamma^{(1)}(\cdot), \mu(\cdot)) + p\widetilde{X}(T, \Gamma^{(2)}(\cdot), \mu(\cdot)).$$

Furthermore, the set **D** is convex; then $\widetilde{X}(t, \Gamma(\cdot), \mu(\cdot)) \in \mathbf{D}$ a.s.. This completes the proof. □

**Lemma A.3** *There exists a constant $c > 0$ such that*

$$\mathbf{E}|\widetilde{X}(T, \Gamma(\cdot), \mu_\alpha(\cdot))|^2 \leq c\left(\|\Gamma(\cdot)\|_{\mathbf{X}}^2 + X_0^2\right) \quad \forall \Gamma(\cdot) \in \mathcal{C}_0^R, \quad \forall \alpha \in \mathcal{T}.$$

*Proof.* For a $\Gamma(\cdot) \in \mathcal{C}_0^R$, let

$$x(t) \triangleq \widetilde{X}(t, \Gamma(\cdot), \mu_\alpha(\cdot)), \quad \pi(t) \triangleq \Gamma(t, [S(\cdot, \mu_\alpha(\cdot)), \mu_\alpha(\cdot)]_{[0,t]}, R_{\mu_\alpha}), \quad \pi(t) = (\pi_1(t), \ldots, \pi_n(t)).$$

By (A.3), it follows that

$$\begin{cases} dx(t) = p(t) \sum_{i=1}^{n} \pi_i(t) \Big( \sum_{j=1}^{n} \sigma_{ij} dw_j(t) + \tilde{a}(t)dt \Big), \\ x(0) = X_0. \end{cases}$$

This is a linear Itô stochastic differential equation, and it is easy to see that the desired estimate is satisfied. This completes the proof. □

**Lemma A.4** *For a given $\alpha \in \mathcal{T}$, the function*

$$\mathbf{E} U(\widetilde{X}(T, \Gamma(\cdot), \mu(\cdot)))$$

*is continuous in $\Gamma(\cdot) \in \mathcal{C}_0^R$.*

*Proof.* Let $\Gamma^{(i)}(\cdot) \in \mathcal{C}_0^R$ and $\widetilde{X}^{(i)}(t) \triangleq \widetilde{X}(t, \Gamma^{(i)}(\cdot), \mu_\alpha(\cdot))$, $i = 1, 2$. By Lemmas A.1 and A.3, it follows that

$$\mathbf{E}|\widetilde{X}^{(1)}(T) - \widetilde{X}^{(2)}(T)|^2 \leq c \|\Gamma^{(1)}(\cdot) - \Gamma^{(2)}(\cdot)\|_{\mathbf{X}}^2,$$

where $c > 0$ is a constant. Then

$$\begin{aligned}
&\left| \mathbf{E} U\left(\widetilde{X}^{(1)}(T)\right) - \mathbf{E} U\left(\widetilde{X}^{(2)}(T)\right) \right| \\
&\leq c_1 \mathbf{E} \left[ (1 + |\widetilde{X}^{(1)}(T)| + |\widetilde{X}^{(2)}(T)|)^{2-q} |\widetilde{X}^{(1)}(T) - \widetilde{X}^{(2)}(T)|^q \right] \\
&\leq c_1 \left[ \mathbf{E} \left(1 + |\widetilde{X}^{(1)}(T)| + |\widetilde{X}^{(2)}(T)|\right)^2 \right]^{1/k'} \left[ \mathbf{E} |\widetilde{X}^{(1)}(T) - \widetilde{X}^{(2)}(T)|^2 \right]^{1/k} \\
&\leq c_2 \left(1 + \|\Gamma^{(1)}(\cdot)\|_{\mathbf{X}}^2 + \|\Gamma^{(2)}(\cdot)\|_{\mathbf{X}}^2\right)^{1/k'} \|\Gamma^{(1)}(\cdot) - \Gamma^{(2)}(\cdot)\|_{\mathbf{X}}^{2/k},
\end{aligned}$$

where $c_i > 0$ are constants, $q$ is as defined in Condition 4.4, $k \triangleq 2/q$, $k' \triangleq k/(k-1) = 2/(2-q)$. This completes the proof. $\square$

Let $S_{\alpha*}(t)$ be defined as the solution of (2.1) with substituting $a_i(t) \equiv M_r(t, \alpha, S(\cdot)|_{[0,t]})$ for all $i$. Clearly, there exists a measurable function $M_\theta(\cdot)$, $M_\theta(t, \cdot) : [0, t] \times \mathcal{T} \times C([0, t]; \mathbf{R}^n) \to \mathbf{R}^m$, such that $\theta_{\mu_\alpha}(t) = M_\theta\left(t, \alpha, S(\cdot)|_{[0,t]}\right)$.

Set $\bar{M}_{\theta*}(t, \alpha) \triangleq M_\theta\left(t, \alpha, S_{\alpha*}(\cdot)|_{[0,t]}\right)$. For an $\alpha \in \mathcal{T}$, set

$$\theta_*(t, \alpha) \triangleq \bar{M}_{\theta*}(t, \alpha).$$

Let

$$z_*(\alpha, T) \triangleq \exp\left(\int_0^T \theta_*(t, \alpha)^\top dw(t) - \frac{1}{2} \int_0^T |\theta_*(t, \alpha)|^2 dt\right).$$

For $\alpha \in \mathcal{T}$, set

$$J'(\Gamma(\cdot), \alpha) \triangleq \mathbf{E} U(\widetilde{X}(T, \Gamma(\cdot), \mu_\alpha(\cdot))).$$

**Lemma A.5** *For a given $\Gamma(\cdot) \in \mathcal{C}_0^R$ and $\alpha^{(2)} \in \mathcal{T}^{(2)}$, the function $J'(\Gamma(\cdot), \alpha)$, where $\alpha = (\alpha^{(1)}, \alpha^{(2)})$, is continuous in $\alpha^{(1)} \in \mathcal{T}^{(1)}$.*

*Proof.* By Condition 4.5(iii), it suffices to consider case $p < 2$ only. Let $\Gamma(\cdot) \in \mathcal{C}_0^R$ and $\alpha^{(2)} \in \mathcal{T}^{(2)}$ be fixed. Let $\alpha_i^{(1)} \in \mathcal{T}^{(1)}$, $i = 1, 2$, $\alpha_i \triangleq (\alpha_i^{(1)}, \alpha^{(2)})$. Set

$$Y_\alpha \triangleq \widetilde{X}(T, \Gamma(\cdot), \mu_\alpha(\cdot)), \quad \alpha \in \mathcal{T}, \qquad Y_* \triangleq \widetilde{X}(T, \Gamma(\cdot), \mu_*(\cdot)),$$

where $\mu_*(t) \triangleq [r(t), 0, \sigma(t)]$. By Girsanov's Theorem (see, e.g., Gihman and Skorohod (1979)), it follows that

$$\begin{aligned}
|\mathbf{E}U(Y_{\alpha_1}) - \mathbf{E}U(Y_{\alpha_2})| &= |\mathbf{E}[z_*(\alpha_1, T) - z_*(\alpha_2, T)]U(Y_*)| \\
&\leq c_1 \mathbf{E}|z_*(\alpha_1, T) - z_*(\alpha_2, T)|(|Y_*|^p + 1) \\
&\leq c_2 \left(\mathbf{E}|z_*(\alpha_1, T) - z_*(\alpha_2, T)|^{q'}\right)^{1/q'} (\mathbf{E}|Y_*|^p + 1)^q)^{1/q} \\
&\leq c_3 \left(\mathbf{E}|z_*(\alpha_1, T) - z_*(\alpha_2, T)|^{q'}\right)^{1/q'} (\mathbf{E}|Y_*|^2 + 1)^q)^{1/q} \\
&\leq c_4 \left(\mathbf{E}|z_*(\alpha_1, T) - z_*(\alpha_2, T)|^{q'}\right)^{1/q'} (\|\Gamma(\cdot)\|_{\mathbf{X}}^2 + 1)^{1/q},
\end{aligned}$$

where $p \in (1, 2)$ is as defined in Conditions 4.4 and 4.5(ii), $q \triangleq 2/p$, $q' \triangleq q/(q-1)$ and $c_i > 0$ are constants.

Furthermore, it is easy to see that for an $\alpha \in A$, we have $z_*(\alpha, T) = y(T)$, where $y(t) = y(t, \alpha)$ is the solution of the equation

$$\begin{cases} dy(t) = y(t) \bar{M}_{\theta*}(t, \alpha)^\top dw(t), \\ y(0) = 1. \end{cases}$$

It is well known that $y(T)$ depends on $\alpha \in \mathcal{T}$ continuously in $L^{q'}(\Omega, \mathcal{F}, \mathbf{P})$ (see, e.g., Krylov (1980, Ch.2)). Hence

$$\mathbf{E}|z_*(\alpha_1, T) - z_*(\alpha_2, T)|^{q'} \to 0 \quad \text{as} \quad \alpha_1 \to \alpha_2.$$

This completes the proof. $\square$

Let $\mathcal{V}$ be the set of all $\sigma$-additive probability measures on $\mathcal{F}_\mathcal{T}$. We consider $\mathcal{V}$ as a subset of $C(\mathcal{T}; \mathbf{R})^*$. (Remind that the set $\mathcal{T}^{(2)}$ is at most countable; in fact, we mean that $C(\mathcal{T}; \mathbf{R})$ is a subspace of the space of all bounded functions mapping $\mathcal{T}$ to $C(\mathbf{R})$ and has the same topology as the space of all bounded functions mapping $\mathcal{T}^{(2)}$ to $C(\mathcal{T}^{(1)}; \mathbf{R})$.) Let $\mathcal{V}$ be equipped with the weak* topology in the sense that

$$\nu_1 \to \nu_2 \quad \Leftrightarrow \quad \int_\mathcal{T} \nu_1(d\alpha) f(\alpha) \to \int_\mathcal{T} \nu_2(d\alpha) f(\alpha) \quad \forall f(\cdot) \in C(\mathcal{T}; \mathbf{R}).$$

**Lemma A.6** *The set $\mathcal{V}$ is compact and convex.*

*Proof.* The convexity is obvious. It remains to show the compactness of the set $\mathcal{V}$. In our case, $\mathcal{T}$ is a compact subset of finite-dimensional Euclidean space. Now we note that the Borel $\sigma$-algebra of subsets of $\mathcal{T}$ coincides with the Baire $\sigma$-algebra (see, e.g., Bauer (1981)). Hence, $\mathcal{V}$ is the set of Baire probability measures. By Theorem IV.1.4 from Warga (1972), it follows that $\mathcal{V}$ is compact. This completes the proof. □

We are now in the position to give a proof of Theorem A.1.

*Proof of Theorem A.1.* For a $\Gamma(\cdot) \in \mathcal{C}_0^R$, we have $J'(\Gamma(\cdot), \cdot) \in C(\mathcal{T}; \mathbf{R})$ and

$$\begin{aligned}\mathbf{E}U(\widetilde{X}(T, \Gamma(\cdot), \mu(\cdot))) &= \int_\mathcal{T} d\nu_\mu(\alpha) \mathbf{E}U(\widetilde{X}(T, \Gamma(\cdot), \mu_\alpha(\cdot))) \\ &= \int_\mathcal{T} d\nu_\mu(\alpha) J'(\Gamma(\cdot), \alpha),\end{aligned}$$

where $\nu_\mu(\cdot)$ is the measure on $\mathcal{T}$ generated by $\Theta$ which corresponds $\mu(\cdot)$. Hence, $\mathbf{E}U(\widetilde{X}(T, \Gamma(\cdot), \mu(\cdot)))$ is uniquely defined by $\nu_\mu$. Let

$$J(\Gamma(\cdot), \nu_\mu) \triangleq \mathbf{E}U(\widetilde{X}(T, \Gamma(\cdot), \mu(\cdot))).$$

By Lemma A.5, $J(\Gamma(\cdot), \nu)$ is linear and continuous in $\nu \in \mathcal{V}$ given $\Gamma(\cdot)$.

To complete the proof, it suffices to show that

$$\sup_{\Gamma(\cdot) \in \mathcal{C}_0^R} \inf_{\nu \in \mathcal{V}} J(\Gamma(\cdot), \nu) = \inf_{\nu \in \mathcal{V}} \sup_{\Gamma(\cdot) \in \mathcal{C}_0^R} J(\Gamma(\cdot), \nu). \tag{A.4}$$

We note that $J(\Gamma(\cdot), \nu) : \mathcal{C}_0^R \times \mathcal{V} \to \mathbf{R}$ is linear in $\nu$. By Lemmas A.1 and A.4-A.5, it follows that $J(\Gamma(\cdot), \nu)$ is either concave or convex in $\Gamma(\cdot)$ and that $J(\Gamma(\cdot), \nu) : \mathcal{C}_0^R \times \mathcal{V} \to \mathbf{R}$ is continuous in $\nu$ for each $\Gamma(\cdot)$ and continuous in $\Gamma(\cdot)$ for each $\nu$. Furthermore, $\mathcal{C}_0^R$ and $\mathcal{V}$ are convex and $\mathcal{V}$ is compact. By the Sion Theorem (see, e.g., Parthasarathy and Ragharan (1971, p.123)), it follows that (A.4), and hence (A.2), are satisfied. This completes the proof of Theorem A.1. □

We are now in the position to give a proof of Theorem 5.1.

## A.3 Proof of Theorem 5.1

Let $\widehat{\alpha} \in \mathcal{T}$ be such that $R_{\widehat{\mu}} = R_{min}$, where $\widehat{\mu}(\cdot) \triangleq \mu_{\widehat{\alpha}}(\cdot)$. By Lemma 5.2(iii)-(iv), it follows that

$$\mathbf{E}U(\widetilde{X}(T, \widehat{\Gamma}_{\widehat{\alpha}}(\cdot, R_{\widehat{\mu}}), \widehat{\mu}(\cdot))) \leq \mathbf{E}U(\widetilde{X}(T, \widehat{\Gamma}_\alpha(\cdot, R_{\mu_\alpha}), \mu_\alpha(\cdot))) \quad \forall \alpha \in \mathcal{T}. \tag{A.5}$$

(If $R_{\widehat{\mu}} = R_{\mu_\alpha}$, then statement (iii) is applicable; if $R_{\widehat{\mu}} < R_{\mu_\alpha}$, then statement (iv) is applicable).

Let $\mu(\cdot) \in \mathcal{A}(\mathcal{T})$ be arbitrary, and let $\nu_\mu(\cdot)$ be the measure on $\mathcal{T}$ generated by $\Theta$, which corresponds to $\mu(\cdot)$. By (A.5), it follows that

$$
\begin{aligned}
\mathbf{E}U(\widetilde{X}(T,\widehat{\Gamma}_{\widehat{\alpha}}(\cdot,R_{\widehat{\mu}}),\widehat{\mu}(\cdot))) &\leq \int_\mathcal{T} d\nu_\mu(\alpha) \mathbf{E}U(\widetilde{X}(T,\widehat{\Gamma}_\alpha(\cdot,R_{\mu_\alpha}),\mu_\alpha(\cdot))) \\
&= \sup_{\Gamma(\cdot) \in \mathcal{C}_0^R} \mathbf{E}U(\widetilde{X}(T,\Gamma(\cdot),\mu(\cdot))) \quad \forall \mu(\cdot) \in \mathcal{A}(\mathcal{T}).
\end{aligned}
\quad (A.6)
$$

By (5.4), (A.6), and Theorem A.1 it follows that the pair $(\widehat{\mu}(\cdot),\widehat{\Gamma}_{\widehat{\alpha}}(\cdot))$ is a saddle point for the problem (3.1)-(3.2). This completes the proof of Theorem 5.1. $\square$

*Proof of Corollary 5.1* follows from Lemma 5.2(iii). $\square$